\documentclass[12pt]{amsart}
\usepackage{amsfonts}
\usepackage{amsmath}
\usepackage{amsxtra}
\usepackage{amssymb,latexsym}
\usepackage[mathcal]{eucal}
\usepackage{amscd}

\newtheorem{theo}{{\bfseries Theorem}}[section]
\newtheorem{prop}[theo]{{\bfseries Proposition}}
\newtheorem{lem}[theo]{{\bfseries Lemma}}
\newtheorem{cor}[theo]{{\bfseries Corollary}}

\theoremstyle{definition}

\newtheorem{df}[theo]{{\bfseries Definition}}

\def \N {\mathbb N}

\def \Z {\mathbb Z}
\def \R {\mathbb R}

\def \B {\mathcal B}

\def \F {\mathcal F}

\def \I {\mathcal I}

\def \M {\mathcal M}

\def \P {\mathcal P}

\def \S {\mathcal S}

\def \00 {\mathbf 0}
\def \a {\alpha }

\def \ep {\epsilon}
\def \g {\gamma}
\def \d {\delta}

\def \proof { {\bfseries Proof:} }

\newcommand{\ga}{\gamma}

\newcommand{\sig}{\sigma}

\usepackage{amssymb,latexsym}
\usepackage[mathcal]{eucal}
\usepackage{amscd}
\usepackage{amsmath}
\usepackage{graphicx}
\usepackage{amsfonts}
\usepackage{amscd}

\numberwithin{equation}{section}

\begin{document}

\title{ Is there a Ramsey-Hindman theorem ?}
\vspace{1cm}


\author{Ethan Akin and Eli Glasner}
\address{Mathematics Department,
 The City College, 137 Street and Convent Avenue,
 New York City, NY 10031, USA}
\email{ethanakin@earthlink.net}

\address{Department of Mathematics,
Tel-Aviv University, Ramat Aviv, Israel}
\email{glasner@math.tau.ac.il}

\date{March, 2015}

\begin{abstract}
We show that there does not exist a joint generalization of the theorems of Ramsey
and Hindman, or more explicitly, that the property of containing a symmetric IP-set is not divisible.
\end{abstract}




\maketitle

%
%

\section{IP and SIP Sets}
\vspace{.5cm}

We will use $\Z, \Z_+, \N$ to stand for the sets of integers, nonnegative integers and positive integers,respectively.

For $F$ a finite subset of  $\Z$, we denote by $\sig_F \in \Z$ the sum of the elements of $F$ with the
convention that $\sig_{\emptyset} = 0$.
Of course, if $F$ is a nonempty subset of $\N$, then $\sig_F \in \N$.

Call a subset $A$ of $\Z$ \emph{symmetric} if  $-A = A$ where $-A = \{ -a : a \in A \}$. For any
subset $A$ of $\Z$ let $A_{\pm} = A  \cup -A $ so that $A_{\pm}$ is the smallest symmetric set
which contains $A $. On the other hand, let $A_+ = A \cap \N$, the positive part of $A$. Note that if
$A$ is symmetric then $A  \ = \ A_{\pm} $ and $A \setminus \{ 0 \} \ = \ (A_+)_{\pm} $.

For subsets $A_1, A_2$ of $\Z$
we let $A_1 + A_2 = \{ a_1 + a_2 : a_1 \in A_1, a_2 \in A_2 \}$ and $A_1 - A_2 = A_1 + (-A_2)$. If
$A_2 = \{ n \}$ we write $A_1 - n$ for $A_1 - A_2$.

 Let $A $ be a nonempty
subset  of  $\Z$.  We set

 \begin{align*}\label{01}
 D(A) & = \{ a_1 - a_2 :  a_1,  a_2 \in A\} = A - A \\
 IP(A)  & = \{ \sig_F : F \ {\text {a finite subset of}} \ A \} \\
 SIP(A) & =  D(IP(A)) =  IP(A) - IP( A) \\
\end{align*}

Clearly, $0 \in D(A)$ and $D(A)$ is symmetric and so the same is true of $SIP(A)$. If $0 \in A$ then $A \subset D(A)$.
In general, $D(A) \cup A \cup -A \ = \ D(A \cup \{ 0 \})$.

In particular, $0 = \sig_{\emptyset}  \in IP(A)$ implies  $IP(A) \subset SIP(A)$.

If $A \subset \N$ then $IP(A) = \{ 0 \} \cup IP(A)_+$ since $0 = \sig_{\emptyset} \in IP(A)$.
$IP(A) \ = \ IP(A \cup \{ 0 \}) \ = \ IP(A \setminus \{ 0 \})$.

\begin{lem}\label{lem101} If $B$ is a nonempty subset of $\N$ then
$$ SIP(B) \quad = \quad  IP(B_{\pm}).$$

If $A$ is a symmetric subset of $\Z$ with $A \setminus \{ 0 \}$ nonempty then
$$ SIP(A_+) \quad = \quad IP(A).$$
\end{lem}

\proof If  $F \subset B_{\pm}$ then
$$\sig_F = \sig_{F \cap B} + \sig_{F \cap -B} = \sig_{F \cap B} - \sig_{(-F) \cap B}$$
Hence,
$IP(B_{\pm}) \subset SIP(B)$

For the reverse inclusion, let $F_1, F_2$ be   finite subsets of $B$.
$$ \sig_{F_1} - \sig_{F_2}   = \sig_{F_1 \cup -F_2} $$
since $F_1$ and $-F_2$ are disjoint.

If $A$ is symmetric and $A \setminus \{ 0 \}$ is nonempty
 then $A_+$ is nonempty and the previous result
applied to $B = A_+$ yields the second  equation since  $(A_{+})_{\pm} =  A \setminus \{ 0 \}$ and
$IP(A) = IP(A \setminus \{ 0 \})$.

$\Box$ \vspace{.5cm}

%
%
%

We say that a subset $B \subset \N$ is
\begin{itemize}
\item  a {\em difference set} if there exists an infinite subset $A$ of $\N$ such that $D(A)_+ \subset B$.
\item an {\em IP set} if there exists an infinite subset $A$ of $\N$ such that $IP(A)_+ \subset B$.
\item an {\em SIP set} if there exists an infinite subset $A$ of $\N$ such that $SIP(A)_+ \subset B$.
\end{itemize}

Since $A \setminus \{ 0 \} = (A_+)_{\pm}$ if $A$ is a symmetric subset of
$\Z$ it follows from Lemma
\ref{lem101} that $B$ is an SIP set iff there
exists an infinite symmetric subset $A$ of $\Z$ such that $IP(A)_+ \subset B$.

We next recall the statements of two -now classical- combinatorial theorems 
(see  \cite{GRS} and \cite{H}):

\begin{theo}\label{theo-ram} [Ramsey]
Let $A$ be an infinite subset of $\N$. If one colors the set $D(A)_+$ in finitely many colors then there
exists an infinite subset $L \subset A$ such that $D(L)_+$ is monochromatic.
\end{theo}

\begin{theo}\label{theo-hind} [Hindman]
Let $A$ be an infinite subset of $\N$. If one colors the set  $IP(A)_+$ in finitely many colors then there
exists an infinite subset  $L \subset \N$ such that $IP(L)_+ \subset IP(A)_+$
and
$IP(L)_+ $
is monochromatic.
\end{theo}
%

In view of these two famous and basic theorems it is 
natural to pose the following question.
Suppose a finite coloring of a set of the form $SIP_+(A) = D(IP(A))_+$ is given, is there an infinite subset
$L \subset \N$ such that $SIP_+(L) \subset SIP_+(A)$ and
$SIP_+(L) $
is monochromatic?

In other words the question is: is there a combined Ramsey-Hindman theorem?

In this
paper
we will show, as expected, that the answer to this question is negative.  We show that it fails
in a strong sense and, in the process, raise some related dynamics questions.
For more details and background see \cite{F} and \cite{A-97}.
We thank Benjy Weiss for his very helpful advice.
\vspace{1cm}

\section{Families of Sets}
\vspace{.5cm}

For an infinite set $Q$ a \emph{family} $\F$ on $Q$ is a collection of subsets of $Q$ which is hereditary upwards.  That is,
$\F \subset \P$, where $\P$ is the power set of $Q$ and
$A \in \F$ and $A \subset B$ implies $B \in \F$. For any collection $\F_1$ of subsets of $Q$, the family
$\F = \{ B : A \subset B $ for some $ A \in \F_1 \}$ is the \emph{family generated by $\F_1$}.

The {\em dual family}
$$ \F_1^*  \ = \ \{ B : B \cap A \not= \emptyset \ \mbox{for all} \ A \in \F_1 \}$$
is indeed a family, and when $\F_1$ is a family we have
 $$ \F_1^* \ = \ \{ B : Q \setminus B \not\in \F_1 \}.$$

A family $\F$ is \emph{proper} when $\F \not= \emptyset$ and $\emptyset \not\in \F$. The dual of a proper family is
proper and $\F^{**} = \F$.

Given families $\F_1, \F_2$ the \emph{join} $\F_1 \cdot \F_2 = \{ A_1 \cap A_2 : A_1 \in \F_1, A_2 \in \F_2 \}$.
By the heredity condition $\F_1 \cup \F_2 \subset \F_1 \cdot \F_2 $.
Clearly, $\F_1 \cdot \F_2$ is proper iff $\F_2 \subset \F_1^*$. We say that two proper families \emph{meet} when
the join is proper.

It is easy to check that for families $\F, \F_1, \F_2$
\begin{itemize}
\item $(\F^*)^* \ = \ \F$.
\item $\F_1 \subset \F_2$ implies $ \F_2^* \subset \F_1^*$. More generally,
$\F \cdot \F_1 \subset \F_2$ implies $\F \cdot  \F_2^* \subset \F_1^*$.
\item $\F^*$ is proper if $\F$ is proper.
\end{itemize}

If $\P_+$ is the collection of all nonempty subsets of $Q$ then $\P_+$ is the largest proper family with dual
$(\P_+)^* = \{ Q \}$, the smallest proper family. The collection $\B$
of all infinite subsets of
$Q$ is a proper family and the dual $\B^*$ is the family of all cofinite subsets of $Q$.

A family is a \emph{filter} when it is proper and closed under finite intersection.  That is, $A_1, A_2 \in \F$ implies
$A_1 \cap A_2 \in \F$. Equivalently, $\F \cdot \F \subset \F$ and so $\F \cdot \F = \F$.  Thus, $\F$ is a filter
iff $\F \cdot  \F^* \subset \F^*$. In particular, if $\F$ is a filter, then $\F \subset \F^*$.
%
%

The dual of a filter is called a \emph{filterdual}. It is sometimes called a \emph{divisible family}. A family is
a filterdual when it satisfies what Furstenberg dubbed the \emph{Ramsey Property}:
\begin{equation}\label{201}
 A_1 \cup A_2 \in \F \qquad \Longrightarrow \qquad A_1 \in \F \quad \mbox{or} \quad A_2 \in \F.
 \end{equation}

A family $\F$ on $Q$ is \emph{full} if it is proper
and $B \in \F$ implies $B \setminus F \in \F$ for any finite  $F \subset Q$ and so
$\F$ is full if $\F \cdot \B^* = \F$. A filter $\F$ is full iff
$\B^* \subset \F$.
In particular, $\B^*$ is the smallest full filter while
$\P_+^*= \{ Q \}$ is the smallest filter.

If a family $\F$ is full then
$$ \F^* \quad = \quad \{ B : B \cap A \ \mbox{is infinite, for all } \ A \in \F \},$$
and $\F^*$ is full.

If $\F$ is a filterdual then, by induction, for all positive integers $k$,
$A_1 \cup \dots \cup A_k \in \F$ implies $A_i \in \F$ for some
$i = 1, \dots, k$. We can interpret this in terms of colorings. If one colors a
set  $A \in \F$ in finitely many colors then there exists $B \in \F$ with $B \subset A$ and $B$ is monochromatic.

Thus, the Ramsey Theorem \ref{theo-ram} implies that the family of difference sets is a filterdual on $\N$ and the
Hindman Theorem \ref{theo-hind} says exactly that the family of IP sets is a filterdual on $\N$.

If $A \subset \N$ is infinite and $K$ is a positive integer then $A \setminus [1,K]$ is infinite and
$IP(A \setminus [1,K])$ is disjoint from $[1,K]$ and is contained in $IP(A)$.  It follows that $\I$ the
family of IP sets is a full family. If $A = \{ a_k : k = 1,2,... \} \subset \N$ with $a_{k+1} > a_k + K$ for all
$k$ then $D(A)_+$ is disjoint from $[1,K]$. Since any infinite set contains such a subsequence it follows that
the family of difference sets is full as well. We will see below that the family $\S$ of SIP sets is also full.

 If $\F$ is a family on $\N$ then $\F$ is \emph{invariant} if $A \in \F$ implies $(A+n)_+ = (A + n) \cap \N \in \F$
 for all $n \in \Z$.
 A proper, invariant family is full since $A \setminus [1,n] = ((A - n)_+ + n)_+$.
 If $\F$ is
 a family of infinite sets, i. e. $\F \subset \B$, then we let
 \begin{equation}\label{202}
 \begin{split}
\g \F \quad = \quad \{ \ (A + n)_+ \ : \ A \in \F, n \in \Z \ \}, \hspace{1cm}\\
\tilde \g \F  \quad = \quad \{ \ A \ : \ (A + n)_+  \in \F,\ \mbox{for all} \  n \in \Z \ \}.
\end{split}
\end{equation}
 That is, $\g \F$ is the smallest invariant family containing $\F$ and $\tilde \g \F$ is a
 largest invariant family contained in $\F$.

 It is easy to see that the dual of an invariant family is invariant from which it follows that for any family $\F$
$$ (\gamma \F)^* \quad = \quad \tilde \gamma (\F^*).$$

Also, one observes that if $\F$ is a filter then $\tilde \g \F$ is an invariant filter contained in $\F$ and
so $\g (\F^*)$ is an invariant filterdual containing the filterdual $\F^*$.

\vspace{1cm}

\section{Dynamics}
\vspace{.5cm}

We call $(X,T)$ a \emph{dynamical system} when $X$ is a compact metric space and $T$ is a homeomorphism on $X$.
We review some well-known facts about such systems.

If $A, B \subset X$ then the \emph{hitting time set} is
$$N(A,B) \ = \ \{ \ n \in \N : T^n(A) \cap B \not= \emptyset \ \}  \ = \ \{ \ n \in \N : A \cap T^{-n}(B) \not= \emptyset \ \}.$$
If $A = \{ x \}$ then we write $N(x,B)$ for $N(A,B)$.
Observe that for $k \in \N$
\begin{equation}\label{301}
\begin{split}
N(A,T^k(B)) \ = \ N(A,B) + k, \hspace{3cm} \\
N(A,T^{-k}(B)) \ \supset \ (N(A,B) - k)_+ \ \supset \ N(A,T^{-k}(B)) \setminus [1,k].
\end{split}
\end{equation}

The system $(X,T)$ is \emph{topologically transitive} if whenever $U, V \subset X$ are nonempty and open,
$N(U,V)$ is nonempty.  In that case, all such $N(U,V)$'s are infinite. A point $x \in X$ is called
a \emph{transitive point} if $N(x,U)$ is nonempty for every open and nonempty $U$ in which case, again,
the $N(x,U)$'s are infinite. We denote by $Trans_T$ the set of transitive points in $X$. The system is
topologically transitive iff $Trans_T$ is nonempty in which case it is a dense $G_{\d}$ subset of $X$.
The system is \emph{minimal} when $Trans_T = X$.

\begin{prop}\label{prop301} If $U, V \subset X$ are nonempty and open and $x$ is a transitive point for $(X,T)$,
then
\begin{equation}\label{302}
N(U,V) \quad = \quad (N(x,V) \ - \ N(x,U))_+. \hspace{2cm}
\end{equation}
\end{prop}

\proof If $n > m$ and $T^n(x) \in V, T^m(x) \in U$ then $T^{n-m}(T^m(x)) = T^n(x)$ implies $n - m \in N(U,V)$.
On the other hand, suppose that $k \in N(U,V)$. Then $U \cap T^{-k}(V)$ is a nonempty open set and so there
exists $m \in \N$ such that $T^m(x) \in U \cap T^{-k}(V)$.  Hence, $T^m(x) \in U$ and $T^n(x) \in V$ with
$n - m = k$.

$\Box$ \vspace{.5cm}

\begin{prop}\label{prop302} Let $U$ be an open set with $x \in U$ where $x$ is a 
recurrent
point for
$(X,T)$. The hitting time set $N(x,U)$ is an IP set.

Assume,in addition, that there exists an involution $J$ on $X$ which maps $T$ to $T^{-1}$ and fixes $x$.  That is,
$J^2 = id_X$, $J \circ T = T^{-1} \circ J$ and $J(x) = x$. In that case, $N(x,U)$ is an SIP set.
\end{prop}

\proof We assume $J$ exists as described above. By intersecting $U$ and $J(U)$ we can assume
that $U$ is $J$ invariant.

 Suppose that $F_N \subset \N$ of cardinality $N$ such that
$SIP(F_N)_+ \subset N(x,U)$. That is, for every $n \in SIP(F_N)_+, \ T^n(x) \in U$. Since $J$ fixes
$x$ and $U$ and maps $T$ to $T^{-1}$ it follows that $T^{-n}(x) \in U$ for all such $n$ as well.
That is, $T^n(x) \in U$ for all $n$ in the symmetric finite set $SIP(F_N)$. Let
$V = \bigcap_{n \in SIP(F_N)}  \  T^{-n}(U)$. By symmetry, $J(V) = V$ and $V$ is a nonempty open set
containing $x$. Since $N(x,V)$ is infinite, there exists $m \in N(x,V)$ which is larger than any element of
$SIP(F_N)$. By construction $T^n(x) \in U$ for all $n \in SIP(F \cup \{ m \})$. It follows that
$F_{N+1} = F_N \cup \{ m \} \subset \N$ has cardinality $N+1$ and $SIP(F \cup \{ m \})_+ \subset N(x,U)$.

Let $F = \bigcup_N \{ F_N \}$. Since we go from $F$ to $IP(F)$ via finite sums, it follows that
 $SIP(F) = \bigcup_N \{ SIP( F_N ) \}$.  Hence, $SIP(F)_+ \subset N(x,U)$.

 The inductive construction for the more general IP result is similar, but easier, as the dance using symmetry is not
 required.

 $\Box$ \vspace{.5cm}

 \begin{cor}\label{cor303} If $(X,T)$ is 
 (recurrent) topologically transitive and $U, V \subset X$ are nonempty and open then
 $N(U,U)$ is an SIP set and $N(U,V)$ is the translation of an SIP set. \end{cor}

 \proof Let $x$ be a transitive point contained in $U$.
 By Proposition \ref{prop301} $N(U,U) = (N(x,U) - N(x,U))_+$ and by Proposition \ref{prop302}
 $N(x,U)$ is an IP set.

 Now let $n \in N(U,V)$ and let $U_0 = U \cap T^{-n}(V)$. $N(U,T^{-n}(V))$ contains the SIP set $N(U_0,U_0)$.
 Since $\S$ is a full family, (\ref{301}) implies that $N(U,V)$ is the translate of an SIP set.

 $\Box$ \vspace{.5cm}

 It is possible to get SIP recurrence under much more general circumstances. 
E.g. as a direct consequence of a theorem of Glasner (\cite[Theorem 1.56]{Gl-03})
we obtain the following result.

\begin{theo}
Let $(X,T)$ be a minimal metric dynamical system, then there is a dense 
$G_\delta$ subset $X_0 \subset X$ such that every $x \in X_0$ is SIP recurrent.
\end{theo}
 
However, Proposition \ref{prop302}  will take care of what we need.
We use it to produce a dynamic example 
which will prove the following:
 
\begin{theo}
The family $\S$ of SIP sets is not a filterdual.
\end{theo}


\begin{proof}
We consider the case where $X$ is the circle $\R/\Z$ and with $a$ a fixed irrational
let $T(x) = x + a$, the \emph{irrational rotation on the circle}. This is a minimal system and so
every point is a transitive point

We can regard the circle as
 $X = [-\frac{1}{2},+\frac{1}{2}]$ with $-\frac{1}{2} = \frac{1}{2}$. The involution $J$ on $X$ given by
 $x \mapsto -x$ fixes $0$ and maps $T$ to $T^{-1}$. Hence, if $U$ is a open set containing $0$ then
 $N(0,U)$ is an SIP set by Proposition \ref{prop302}. If $b \in X$ then the translation $x \mapsto x + b$
 commutes with $T$ and maps $0$ to $b$.  It follows that if $U$ is an open set which contains $b$ then $N(b,U)$ is an
 SIP set.

Let $U = (\frac{1}{8},\frac{1}{8}), U_+ = [0,\frac{1}{8}), U_- = (-\frac{1}{8},0]$.
 The SIP set $N(0,U) $ is the union $ N(0,U_+) \cup N(0,U_-)$ and
 we will show that neither $N(0,U_+)$ nor $N(0,U_-)$ is an SIP set. Replacing $a$ by $-a$ interchanges the two sets and
 so it suffices to focus on $N(0,U_+)$.
 We 
 have to show that there is no infinite subset $A$ of $\N$
such that 
$SIP(A)_+ \subset N_T(0,U_+)$
\vspace{.25cm}

 Assume
 such
 $A$ exists. Let
 $$  M = \sup \ \{ T^t(0) = ta : t \in SIP(A)_+ \}. $$
 Thus, $0 < M \leq \frac{1}{8}$.
 Given any $\ep > 0$ there is a finite subset $F \subset A_{\pm}$ with $0 < \sig_F$ and
 such that $M - \ep < T^{\sig_F}(0) = \sig_F a \leq M \leq \frac{1}{8}$. Since
 $A$ is infinite, there exists $t \in A $ larger than  all the elements of $ SIP(F)_+$ and so with $t >  \sig_F $.
 Thus, $t - \sig_F , t, t + \sig_F \in SIP(A)_+ \subset N_T(0,U_+)$. Thus, $0 < (t - \sig_F)a \leq \frac{1}{8}$.
 Since $2 \sig_F a \leq  2M \leq \frac{1}{4}$, we have $(t + \sig_F )a = (t - \sig_F)a + 2 \sig_F a > 2(M - \ep)$.
 If $\ep$ is chosen less than $\frac{M}{2}$ then  $t + \sig_F \in SIP(A)_+$
 with $(t + \sig_F )a > M$.  This contradicts
 the definition of $M$.
 %
 \end{proof}

$\Box$ \vspace{.5cm}

However, the dynamics suggests a further conjecture.

A dynamical system $(X,T)$ is called $\F$ topologically transitive for a proper family $\F$ of subsets of $\N$
if for all $U, V \subset X$ open and nonempty $N(U,V) \in \F$. From (\ref{301}) it follows that every
translate of $N(U,V)$ is also in $\F$ and so $N(U,V) \in \tilde \g \F$.  That is, an $\F$ topologically
transitive family is automatically a $\tilde \g \F$ topologically transitive family.

A system $(X,T)$ is called \emph{mild mixing} if it is $\S^*$ topologically transitive. Glasner and Weiss
\cite[Theorem 4.11, page 614]{GW} (and also, independently, Huang and Ye \cite{HY}) prove the following.

\begin{theo}\label{prop-mild-mix} $(X,T)$ is mild mixing iff for every topologically transitive system
$(Y,S)$ the product system $(X \times Y, T \times S)$ is topologically transitive. \end{theo}

\proof Suppose $U_1, V_1 \subset X$ and $U_2, V_2 \subset Y$ are open and nonempty.  Fix $n \in N(U_2, V_2)$ so that
$U_3 = U_2 \cap S^{-n}(V_2) \subset Y$ is open and nonempty.  By Corollary \ref{cor303}
$N(U_2,S^{-n}(V_2)) \supset N(U_3,U_3)$ is an SIP set. Because $(X,T)$ is mild mixing
$N(U_1,T^{-n}(V_1))$ is an $SIP^*$ set.  Because $\S$ is a full family the intersection is infinite.
The intersection is $N(U_1 \times U_2, (T \times S)^{-n}(V_1 \times V_2))$ and so by
(\ref{301}), $N(U_1 \times U_2, V_1 \times V_2)$ is infinite.  Thus, the product is topologically transitive.

If $(X,T)$ is not mild mixing then there exist $U,V \subset X$ open and nonempty and an SIP set $A \subset \N$
such that $N(U,V) \cap A = \emptyset$.  The result then follows a construction of 
Glasner-Weiss
which shows

\begin{theo}\label{gw-const} If $A \subset \N$ is an SIP set then there exists a topologically transitive system
$(Y,S)$ and $G \subset Y$ a nonempty open set such that $N(G,G) \subset A$. \end{theo}

$\Box$ \vspace{.5cm}

\begin{cor}\label{cor-mild-mix} The product of any collection of mild mixing systems is mild mixing. \end{cor}

\proof  If $T_1$ and $T_2$ are mild mixing homeomorphisms and $S$ is topologically transitive, then
$T_2 \times S$ is transitive and so $T_1 \times T_2 \times S$ is transitive. Hence,$T_1 \times T_2$ is mild mixing.

By induction a finite product of mild mixing systems is mild mixing.

An infinite product times $S$ is the inverse limit of finite products times $S$ and the inverse limit of
transitive systems is transitive.  It follows that the infinite product is mild mixing.

$\Box$ \vspace{.5cm}

Let $\M$ be the family on $\N$ generated by $\{ N(U,V) : (X,T) $ mild mixing and $U, V \subset X$ open and nonempty $ \}$.
From Corollary \ref{cor-mild-mix} it follows that $\M$ is a filter. Because $\S^*$ transitivity implies
$\tilde \g (\S^*)$ transitivity it follows that $\M \subset \tilde \g (\S^*)$.

By the Hindman Theorem, $\I$ the family of IP sets is a filterdual and so $\I^*$ is a filter. It then follows
that $\tilde \g (\I^*)
= (\ga (\I))^*$ is a filter.

We know from the above example that $\S^*$ is not a filter, but it might still be true that $\tilde \g (\S^*)$ is
a filter.  This would be true if $\M = \tilde \g (\S^*)$. In that case, $\g \S$ would be a filterdual.  In the example itself,
$N(0,U_+)$ is not an SIP set, but if $T^k(0) \in (-\frac{1}{8},0)$ then $0$ is in the interior of $T^k(U_+)$ and so
$N(0,T^k(U_+))$ is an SIP set.  From (\ref{301}) it follows that $N(0,U_+)$ is the translate of an SIP set.

It is to this question that we now turn.  As we will
see this
conjecture fails as well.
\vspace{1cm}


\section{SIP Sets and Their Refinements}
\vspace{.5cm}

Let $e \in \N$ and $b = 2e +1$ so that $b$ is an odd number greater than $1$. Define $\a_b : \N \to \N$ by
$\a_b(n) = b^{n-1}$. The $b$ expansion of an integer $t$ is the sum $\Sigma_{n \in \N} \ \ep_n \a_b(n) \ = \ t$
such that:
\begin{itemize}
\item  $|\ep_n| \leq e$ for all $n \in \N$.
\item  $\ep_n = 0$ for all but finitely many $n$.
\end{itemize}

\begin{prop}\label{prop401}  Every integer in $\Z$ has a unique $b$ expansion.

\end{prop}

\proof  By the Euclidean Algorithm every integer $t$ can be expressed uniquely as $\ep + bs$ with $|\ep| \leq e$.
It follows by induction that every integer $t$ with $|t| <\frac{1}{2}(b^k - 1)$ has an expansion with
$\ep_n = 0$ for $n \geq k$. There are $b^k$ such integers and the same number of expansions. So by the pigeonhole
principle the expansions are unique.

$\Box$ \vspace{.5cm}

We will only need the $b = 3$ expansions with $e = 1$ so that each $\ep_n = -1, +1$ or $0$. We will write
$\a$ for $\a_3$ so that $\a(n) = 3^{n-1}$. From Proposition \ref{prop401} we obviously have
 $$\Z \quad = \quad SIP(\a(\N)).$$

The \emph{length} $r(t)$ of $t$ is the number of nonzero $\ep_i$'s in the expansion of $t$.
With $r = r(t)$ we let $j_1(t),...,j_r(t)$
be the corresponding indices written in increasing order.  That is,
\begin{itemize}
\item $j_1(t) < \dots < j_r(t)$  and $\ep_{j_i(t)} \ = \ \pm 1$ for $i = 1, \dots, r = r(t)$.
\item $t \quad = \quad \Sigma_{i = 1}^{r} \ \ep_{j_i(t)} \a(j_i(t))$.
\end{itemize}
 We call this representation the \emph{reduced expansion} and $j_1(t), \dots , j_r(t)$ the \emph{indices} of $t$.

 Notice that $0$ has length $0$ and equals the empty sum.

 Because $3^{n+1} > 1 + 3 + \dots + 3^n$ it follows that
\begin{equation}\label{401}
 t > 0 \quad \Leftrightarrow \quad \ep_{j_r(t)} = 1.
 \end{equation}

\begin{df}\label{df402}  Assume that
 $j_1(t),\dots,j_{r(t)}(t)$ and $j_1(s), \dots, j_{r(s)}(s)$ are the indices of the reduced expansions for $t, s \in \Z$.

(a) Call $t$ of \emph{positive type} (or of \emph{negative type}) if
$\ep_{j_1(t)}\ep_{j_{r}(t)}$ is positive (resp. is
negative). So $t$ is of positive type if coefficients of its first and last indices have the same sign.
By convention we will say that
$0$ is of positive type.

(b) We will write $t \succ s$ if $j_1(t) > j_{r(s)}(s)$, that is, the indices for $t$ are larger
than all of the indices of $s$.  We will say that $t$ is \emph{beyond} $s$ when $t \succ s$.
\end{df}
\vspace{.5cm}

If $t > 0$ then
$\ep_{j_r(t)} = 1$, and so $t$ is of positive type (or negative type) if
$\ep_{j_1(t)}$ is positive (resp. $\ep_{j_1(t)}$ is negative).

Notice that
 if $j_{r(s)}(s) = n + 1$, then
\begin{equation}\label{402}
t \succ s \qquad \Longleftrightarrow \qquad t \succ  3^n \qquad \Longleftrightarrow \qquad t \ \equiv \ 0 \quad
(mod \ 3^{n+1}),
\end{equation}

Now we turn to SIP sets.

\begin{df}\label{df403} (a) We call a strictly increasing function $k : \N \to \N$ a \emph{+function}.

(b) If $k_1$ and $k_2$ are +functions we say that $k_2$ \emph{directly refines} $k_1$ if $k_2(\N) \subset k_1(\N)$.
We say that $k_2$ \emph{refines} $k_1$ if $IP(k_2(\N)) \subset IP(k_1(\N))$.

\end{df}

Clearly, direct refinement implies refinement and each relation is transitive.

For an infinite subset $A \subset \N$ there is a unique +function $k_A$ such that $k_A(\N) = A$, i. e. the function which
counts the elements of $A$ in increasing order.

\begin{lem}\label{lem404} If $k$ is a +function, then for any $N \in \N$ there exists a +function $k_1$ such that
\begin{itemize}
\item  $k_1$ refines $k$.
\item $k_1(n) \succ 3^{N-1}$ for all $n \in \N$.
\end{itemize}
\end{lem}

\proof By induction we can assume that $k(n) \succ 3^{N-2}$ for all $n \in \N$ (the condition is vacuous when $N = 1$).
This means that $j_1(k(n)) \geq N$ for all $n  \in \N$.

Case 1: There is an infinite set $A \subset \N$ such that $j_1(k(n)) > N$ for all $n \in A$. Let $k_1$ be the
+function with $k_1(\N) = k(A)$. Then $k_1$ is a direct refinement of $k$ and $j_1(k_1(n)) \geq N+1$ for all $n \in \N$,
i. e. $k_1(n) \succ 3^{N - 1}$ for all $n$.

Case 2: There is an infinite set $A \subset \N$ such that $j_1(k(n)) = N$ and
$\ep_{j_1(k(n))} = \d = -1$ for all $n \in A$, or
there is an infinite set $A \subset \N$ such that $j_1(k(n)) = N$ and
$\ep_{j_1(k(n))} = \d = +1$ for all $n \in A$.

Let $\tilde k$ be the +function with $\tilde k(\N) = k(A)$, a direct refinement of $k$ and
$\tilde k(n) \equiv \d 3^{N-1} \ (mod \ 3^N)$ for all $n \in \N$. 

Now define
$$ k_1(n) = \tilde k(3n - 2) + \tilde k(3n -1) + \tilde k(3n).$$
Clearly, $IP(k_1(\N)) \subset IP(\tilde k(\N)) \subset IP(k(\N))$ and so $k_1$ refines $k$. For all $n$,
$k_1(n) \equiv 3\d 3^{N-1} \equiv 0$
$\ (mod \ 3^N)$. Thus, $k_1(n) \succ 3^{N - 1}$ for all $n$.

$\Box$ \vspace{.5cm}

\begin{theo}\label{theo405} If $A \subset \N$ is a translate of an SIP set then
there exists $t_0 \in A$ and +function $k$ such that
\begin{itemize}
\item[(i)] $k(1) \succ t_0$.
\item[(ii)] $k(n+1) \succ k(n)$ for all $n \in \N$.
\item[(iii)] Either $k(n)$ is of positive type for all $n \in \N$, or else
$k(n)$ is of negative type for all $n \in \N$.
\item[(iv)] $t_0 + SIP(k(\N))_+ \ = \ (t_0 + SIP(k(\N)))_+ \subset A$.
\end{itemize}
\end{theo}

\proof There exists $u \in \Z$ and a +function  $k_0$ such that $(SIP(k_0(\N)) + u)_+ \subset A$.

For sufficiently large $N_0$,
$t_0 = u + \Sigma_{n=1}^{N_0} k_0(n)  > 0$ and so lies in $A$.

Let $k^+_0$ be the direct refinement of $k_0$ with $k^+_0(\N) = k_0([N_0+1,\infty))$.
Hence,
\begin{equation}\label{403}
(t_0 + SIP(k^+_0(\N)))_+ \ \subset \ A.
\end{equation}

Now we repeatedly apply Lemma \ref{lem404}.

Let $N_1 > j_r(t_0)$.

Choose $k_1$ a +function which refines $k^+_0$ and with $k_1(n) \succ 3^{N_1}$ for all $n \in \N$. In particular,
$k_1(1) \succ t_0$. So from (\ref{403}) we have
\begin{equation}\label{404}
(t_0 + SIP(k_1(\N)))_+ \ \subset \ A.
\end{equation}
Let $k_1^+$ be the direct refinement of $k_1$ with $k_1^+(\N) = k_1([2,\infty))$.

Let $N_2 > j_r(k_1(1))$ and choose $k_2$ a +function which refines $k^+_1$ and with
$k_2(n) \succ 3^{N_2}$ for all $n \in \N$. In particular,
$k_2(1) \succ k_1(1)$.  Furthermore,
\begin{equation}\label{405}
IP[\{ k_1(1) \} \cup IP(k_2(\N))] \subset IP[\{ k_1(1) \} \cup IP(k_1^+(\N))]  = IP(k_1(\N)).
\end{equation}

Inductively, let $k_q^+$ be the direct refinement of $k_q$ with $k_q^+(\N) = k_q([2,\infty)$ and let $N_{q+1} > j_r(k_q(1))$.
Choose $k_{q+1}$ a refinement of $k_q^+$ with $k_{q+1}(n) \succ 3^{N_{q+1}} $
for all $n \in \N$.  Hence, $k_{q+1}(1) \succ k_q(1)$ and
\begin{equation}\label{406}
\begin{split}
IP[\{ k_q(1) \} \cup IP(k_{q+1}(\N))] \subset IP[\{ k_q(1) \} \cup IP(k_q^+(\N))]  = IP(k_q(\N)), \\
IP[\{k_1(1), \dots,  k_q(1) \} \cup IP(k_{q+1}(\N))] \subset  IP(k_1(\N)).
\end{split}
\end{equation}

Now define $\tilde k(n) = k_n(1)$ for $n \in \N$. Either $\tilde k(n)$ is of positive type infinitely often or of
negative type infinitely often (or both). So we can choose a direct refinement $k$ of $\tilde k$ so that, (i), (ii)
and (iii) hold.  In addition,
$$IP(k(\N)) \subset IP(\tilde k(\N)) \subset IP(k_1(\N)).$$

Clearly, $k(n+1) \succ k(n) \succ t_0$ and from
(\ref{406}) it follows that $IP(\tilde k(\N)) \subset IP(k_1(\N))$ and so from (\ref{404})
$[t_0 + SIP(k(\N))]_+ \subset A$.

Since $k(n) \succ t_0$ for all $n$ it follows from (\ref{402}) that $t \succ t_0$ for all $t \in SIP(k(\N))$.
Hence, if $t \in SIP(k(\N))$ is negative then $t_0 + t$ is negative.  Thus, $[t_0 + SIP(k(\N))]_+ =  t_0 + [SIP(k(\N))]_+$,
completing the proof of (iv).

$\Box$ \vspace{.5cm}

{\bfseries Remark:}  Since  $N_1$ can be chosen arbitrarily large it follows that  
$SIP(k(\N))_+ $ and hence $t_0 + SIP(k(\N))_+ $ can be chosen disjoint from an arbitrary finite subset of $\N$.
This shows that $\S$ and $\gamma \S$ are full families.
\vspace{.5cm}

 For two distinct numbers $n, m \in \Z \setminus \{ 0 \}$ define
 \begin{equation}\label{407}
 \d(n,m) \quad = \quad \begin{cases} 0 \quad \mbox{if} \ nm > 0, \\
 1 \quad \mbox{if} \ nm < 0.\end{cases}
 \end{equation}

 Now we define the \emph{sign change count} to be the function $z : \N \to \Z_+$
 so that if $  t \in \N$ has reduced expansion with indices $j_1(t),\dots,j_{r(t)}(t)$ then
 \begin{equation}\label{408}
 z(t) \quad = \quad \Sigma_{i = 1}^{r(t)-1} \ \d(\ep_{j_i},\ep_{j_{i+1}}).
 \end{equation}
 In particular, if the length is  one then the sum is empty and so
 $z(3^{n-1}) = 0$ for all $n \in \N$.

 For a positive integer $K$ let $\pi_K : \Z \to \Z/K \Z$ be the quotient map mod $K$.

 \begin{theo}\label{theo406}
 If $A \subset  \N$ is a translate of an SIP set
 then for every odd number $K$, $\pi_K \circ z : A \to \Z/K \Z$ is surjective. \end{theo}

 \proof Fix $K$. Since it is odd, $2$ and $-2$ generate the cyclic group $\Z/K \Z$.

 By Theorem \ref{theo405} we can choose $t_0 \in A$ and a +function $k$ which satisfies the
 four conditions of the theorem.

Let $s_0 = t_0 + \Sigma_{n=1}^{2K+1} \ k(n)$. Since $k(n+1) \succ k(n) \succ t_0$ for all $n$ we can regard the
sequence $\{ k(n) : n \in \N \}$ as a sequence of disjoint ascending blocks in $IP(\a)$.

Since each $k(n)$ is positive, each
$\ep_{j_{r_n}(k(n))}$ is positive, where $r_n = r(k(n))$.
For
$i = 1,...,K-1$ let
$$ s_i = t_0 + \Sigma_{n=1}^{2i} (-1)^{n+1} k(n) + \Sigma_{n = 2i +1}^{2K +1} k(n).$$
That is, moving from $s_{i-1}$ to $s_{i}$ we reverse the sign of the block
$k(2i)$
keeping the remaining blocks fixed.
Clearly $s_i \in A$ for $i = 1,...,K-1$.
\vspace{.25cm}

Case 1- Every $k(n)$  is of positive type.  Each
$\ep_{j_1(k(n))}$ is positive.

Moving from $s_{i-1}$ to $s_{i}$ increases $z$ by exactly $2$
because the $++$ transition from $j_{r_{2i-1}}(k(2i-1))$ to $j_{1}(k(2i))$ is replaced by a $+-$ transition and
the $++$ transition from $j_{r_{2i}}(k(2i))$ to $j_{1}(k(2i+1))$ is replaced by a $-+$ transition.
Thus,
$\pi_K(z(s_i)) = \pi_K(z(s_{i-1})) + 2 \
(mod \ K)$. Since $2$ generates the cyclic group, $\pi_K \circ z$ is surjective.
\vspace{.25cm}

Case 2- Every $k(n)$ for $n > 1$ is of negative type. Each
$\ep_{j_{1}(k(n))}$
is negative.
This time the $+-$ transition from $j_{r_{2i-1}}(k(2i-1))$ to $j_{1}(k(2i))$ is replaced by a $++$ transition
and the $+-$  transition from $j_{r_{2i}}(k(2i))$ to $j_{1}(k(2i+1))$  is replaced by a $--$ transition.
Thus, in this case,
$\pi_K(z(s_i)) = \pi_K(z(s_{i-1})) - 2 \ (mod \ K)$.
Again $\pi_K \circ z$ is surjective.

$\Box$ \vspace{.5cm}

We can now deduce the following:

\begin{theo}
\label{407a}
If $A$ is any SIP subset of $\N$ (including $\N$ itself), then
$A$ can be partitioned by two sets neither of which contains a translate of an SIP set.
Thus, the family of translated SIP sets in $\N$ is not a filterdual.
 \end{theo}

\proof With $K = 3$, the sign count map $z : \N \to \Z/3 \Z$ determines a coloring of $\N$ and in any
translated SIP set contained in $IP(k)_+$ all three colors occur.

In particular, let
\begin{equation}\label{409} \begin{split}
A_0 \quad = \quad \{ \ t \in \N \ : \ z(t) \equiv 0 \
(mod \ 3) \ \}, \hspace{3cm}\\
 A_1 \quad = \quad \N \setminus A_0 \quad = \quad  \{ \ t \in \N \ : \ z(t) \not\equiv 0 \
 (mod \ 3) \ \}.\hspace{1cm}
\end{split}
\end{equation}
Neither $A_0$ nor $A_1$ contains a translate of an SIP set. It follows that both $A_0$ and $A_1$ are elements of
the dual $(\gamma \S)^* = \tilde \gamma (\S^*)$.

$\Box$ \vspace{1cm}

In general the congruence classes $mod \ K$ of $z(t)$ (for $K$ odd) define a decomposition of $\N$ into $K$
elements, each a member of  $\tilde \gamma (\S^*)  \subset \S^*$. Thus,
 $\S^*$ and $\tilde \g (\S^*)$ fail to be filters in a very strong way.

\vspace{1 cm}

\section{Dynamics again}
\vspace{.5cm}

We defined the family $\M$ generated by the sets $N(U,V)$ with $(X,T)$ mild mixing and $U,V \subset X$ open
and nonempty. We saw that $\M$ is an invariant filter contained in $\tilde \g (\S^*)$.  Now that we know that
the latter is not a filter, we see that the inclusion is proper. Can we find another possible description of the sets
in $\M$?

For a proper family $\F$ on an infinite set $Q$ we define the \emph{sharp dual} $\F^{\#}$ by
\begin{equation}\label{501}
 \F^{\#} \quad = \quad \{ \ A \subset Q \ : \ A \cap B \in \F \quad \mbox{for all} \ B \in \F \ \}.
\end{equation}

\begin{prop}\label{prop501} Let $\F$ be a proper family on an infinite set $Q$.

(a) $\F^{\#}$ is a filter contained in $\F \cap (\F^*)$. It is full if $\F$ is full.

(b) $\F^{\#} = (\F^*)^{\#}$.

(c) If $\F$ is a filter, then $\F^{\#} = \F$. In particular, $ (\F^{\#})^{\#} \ = \ \F^{\#}$.

(d) If $\F$ is a filterdual then $ \F^{\#} = \F^*$.
\end{prop}

\proof (a) Since $\F$ is proper, $\emptyset \not\in \F$ and so
 $A \in \F^{\#}$ implies $A \cap B \not= \emptyset$ for all $B \in \F$. Thus,
$A \in \F^*$.  Also, $Q \in \F$ and so $A = A \cap Q \in \F$.  Thus, $A \in \F$.

If $A_1, A_2 \in \F^{\#}$ and $B \in \F$ then $(A_1 \cap A_2) \cap B = A_1 \cap (A_2 \cap B) \in \F$.
Thus, $A_1 \cap A_2 \in \F$.

If $A$ is a cofinite set and $B \in F$ then $A \cap B \in \F$ since $\F$ is full.  Thus, $A \in  \F^{\#}$.
That is, $\B^* \subset \F^{\#}$. Since the latter is a filter, it is full.

(b) If $A \in \F^{\#}, B_1 \in \F^*, B \in \F $, then $(A \cap B_1) \cap B = (A \cap B) \cap B_1 \not= \emptyset$.
Since $B$ was arbitrary, $A \cap B_1 \in \F^*$. Since $B_1$ was arbitrary, $A \in  (\F^*)^{\#}$. The reverse
inclusion follows from $(\F^*)^* = \F$.

(c) If $\F$ is a filter, then $\F \cdot \F = \F$ and so $\F \subset \F^{\#}$.  From (a) it follows that
$\F^{\#} \subset \F$.

(d) This follows from (b) and (c).

$\Box$ \vspace{.5cm}

\begin{theo}\label{theo502} $\M \subset (\g \S)^{\#} =  (\tilde \g (\S^*))^{\#} $. \end{theo}

\proof The equation follows from Proposition \ref{prop501} (b).

Now let $(X,T)$ be mild mixing and $U,V \subset X$ be open and nonempty. Let $A \subset \N$ be a translation of
an SIP set.  We show that $N(U,V) \cap A$ is the translation of an SIP set.

By the Glasner Weiss construction Theorem \ref{gw-const} 
there exists a topologically transitive
system, $(Y,S), \ G \subset Y$ open and nonempty and $n \in \Z$ so that $N(G,S^{-n}(G)) \setminus [1,|n|]$
is contained in $A$. It follows that $N(U,V) \cap A $ contains $N(U,V) \cap N(G,T^{-n}(G) \setminus [1,|n|]$.
Because $(X \times Y,T \times S)$ is topologically transitive,
$N(U,V) \cap N(G,T^{-n}(G)) = N(U \times G, V \times T^{-n}(G))$ is the translation of an SIP set by Corollary
\ref{cor303}. As $\g \S$ is a full family, it follows that $N(U,V) \cap A$ is in $\g \S$.

$\Box$ \vspace{1cm}

 Our final conjecture is that $\M = (\g \S)^{\#}$.



\bibliographystyle{amsplain}

\begin{thebibliography}{10}

\bibitem{A-97}
E. Akin,
{\em Recurrence in topological dynamical systems:
Furstenberg families and Ellis actions}, Plenum Press, New York,
1997.



\bibitem{F}
H. Furstenberg,
{\em Recurrence in ergodic theory and combinatorial number theory},
Princeton university press,  Princeton, N.J., 1981.




\bibitem{Gl-03}
E. Glasner,
{\it Ergodic theory via joinings}, Math. Surveys and
Monographs, AMS, {\bf 101}, 2003.

%
\bibitem{GW}
E. Glasner and B. Weiss,
{\em On the interplay between measurable and topological dynamics},
Handbook of dynamical systems. Vol. {\bf 1B}, 597--648, Elsevier B. V., Amsterdam, 2006.

\bibitem{GRS}
R. Graham, B. L. Rothchild and J. H. Spencer,
{\em Ramsey theory},
John Wiley \& sons, 1980.

\bibitem{H}
N. Hindman, 
{\em  Finite sums from sequences within cells of a partition of N},
 J. Combinatorial Theory, Ser. A {\bf 17}, (1974), 1--11. 


\bibitem{HY}
W.  Huang and X. Ye,
{\em Topological complexity, return times and weak disjointness},
Ergodic Theory Dynam. Systems {\bf 24}, (2004), 825--846.

%

\end{thebibliography}

 \end{document}